\input amstex\documentstyle{amsppt}  
\pagewidth{12.5cm}\pageheight{19cm}\magnification\magstep1
\topmatter
\title Elliptic Weyl group elements and unipotent isometries with $p=2$
\endtitle
\author George Lusztig and Ting Xue\endauthor
\address{Department of Mathematics, M.I.T., Cambridge, MA 02139;
Department of Mathematics, Northwestern University, Evanston, IL 60208}\endaddress
\thanks{G.L. supported in part by the National Science Foundation}\endthanks
\endtopmatter   
\document

\define\uuG{\un{\un G}}

\define\uWW{\un\WW}

\define\dsv{\dashv}

\define\pe{\perp}

\define\sqc{\sqcup}

\define\qua{\quad}

\define\op{\oplus}

\define\part{\partial}
\define\em{\emptyset}

\define\m{\mapsto}
\define\do{\dots}

\define\sub{\subset}    

\define\T{\times}

\define\nl{\newline}
\redefine\i{^{-1}}

\define\un{\underline}

\define\g{\gamma}
\redefine\d{\delta}
\define\e{\epsilon}

\define\s{\sigma}

\redefine\L{\Lambda}
\define\Ph{\Phi}

\define\kk{\bold k}

\define\nn{\bold n}

\define\NN{\bold N}

\define\WW{\bold W}
\define\ZZ{\bold Z}

\define\cf{\Cal F}

\define\ct{\Cal T}

\define\fB{\frak B}

\define\fS{\frak S}

\define\WE{L1}
\define\WEH{L2}
\define\SPA{S}

\head Introduction \endhead
\subhead 0.1\endsubhead
Let $G$ be a connected reductive algebraic group over an algebraically closed field $\kk$ of characteristic 
$p\ge0$. Let $\uuG$ be the set of unipotent conjugacy 
classes in $G$. Let $\uWW$ be the set of conjugacy classes in the Weyl group $\WW$ of $G$. 
For $w\in\WW$ and $\g\in\uuG$ let $\fB^\g_w$ be the variety of all pairs $(g,B)$ where $g\in\g$ and $B$ is a
Borel subgroup of $G$ such that $B$ and $gBg\i$ are in relative position $w$. For $C\in\uWW$ and
$\g\in\uuG$ we write $C\dsv\g$ when for some (or equivalently any) element $w$ of minimal length in $C$ we have
$\fB^\g_w\ne\em$. In \cite{\WE, 4.5} a natural surjective map $\Ph:\uWW@>>>\uuG$ was defined. When $p$ is not a 
bad prime for $G$, the map $\Ph$ can be characterized in terms of the relation $C\dsv\g$ as follows (see 
\cite{\WE, 0.4}):

(a) {\it If $C\in\uWW$ then $\Ph(C)$ is the unique unipotent class of $G$ such that $C\dsv\Ph(C)$ and such that 
if $\g'\in\uuG$ satisfies $C\dsv\g'$ then $\Ph(C)$ is contained in the closure of $\g'$.}
\nl
If $p$ is a bad prime for $G$ then the definition of the map $\Ph$ given in \cite{\WE} is less direct; one first
defines $\Ph$ on elliptic conjugacy classes by making use of the analogous map in characteristic $0$ and then one
extends the map in a standard way to nonelliptic classes. It would be desirable to establish property (a) also 
in bad characteristic. To do this it is enough to establish (a) in the case where $C$ is elliptic (see the
argument in \cite{\WE, 1.1}.) One can also easily reduce the general case to the case where $G$ is almost simple;
moreover it is enough to consider a single $G$ in each isogeny class. The fact that (a) holds for $C$ elliptic
with $G$ almost simple of exceptional type (with $p$ a bad prime) was pointed out in \cite{\WEH, 4.8(a)}.
It remains then to establish (a) for $C$ elliptic in the case where $G$ is a symplectic or special 
orthogonal group and $p=2$. This is achieved in the present paper. In fact, Theorem 1.3 establishes (a) with $C$ 
elliptic in the case where $G$ is $Sp_{2n}(\kk)$ or $SO_{2n}(\kk)$ ($p=2$); then (a) for $G=SO_{2n+1}(\kk)$ 
($p=2$) follows from the analogous result for $Sp_{2n}(\kk)$ using the exceptional isogeny 
$SO_{2n+1}(\kk)@>>>Sp_{2n}(\kk)$. Thus the results of this paper establish (a) for any $G$ without restriction on
$p$.

\subhead 0.2\endsubhead
If $w\in\WW$ and $\g\in\uuG$ then $G_{ad}$ (the adjoint group of $G$) acts on $\fB^\g_w$ by 
$x:(g,B)\m(xgx\i,xBx\i)$.
Let $C\in\uWW$ be elliptic. Let $\g=\Ph(C)$. The following result is proved in \cite{\WEH, 0.2}.

(a) {\it For any $w\in C$ of minimal length, $\fB^\g_w$ is a single $G_{ad}$-orbit.}
\nl
The following converse of (a) appeared in \cite{\WEH, 3.3(a)} in the case where $p$ is not a bad prime for $G$ and
in the case where $G$ is almost simple of exceptional type and $p$ is a bad prime for $G$ (see also
\cite{\WE, 5.8(c)}):

(b) {\it Let $\g'\in\uuG$. If $C\dsv\g'$ and $\g'\ne\Ph(C)$ then for any $w\in C$ of minimal length,
$\fB^\g_w$ is a union of infinitely many $G_{ad}$-orbits.}
\nl
Using 0.1(a) we see as in the proof of \cite{\WE, 5.8(b)} that (b) holds for any $G$ without restriction on $p$.
Namely, from \cite{\WE, 5.7(ii)} we see that $\fB^{\g'}_w$ has pure dimension equal to $\dim\g'+l(w)$ where $l(w)$
is the length of $w$ and $\fB^\g_w$ has pure dimension equal to $\dim\g+l(w)$. Also by \cite{\WE, 5.2} the action
of $G_{ad}$ on $\fB^{\g'}_w$ or $\fB^\g_w$ has finite isotropy groups. Thus, 
$\dim\fB^\g_w=\dim G_{ad}$ (see (a)) and to prove (b) it is enough to show that $\dim\fB^{\g'}_w>\dim G_{ad}$ or 
equivalently that $\dim\g'+l(w)>\dim\g+l(w)$ or that $\dim\g'>\dim\g$. But from 0.1(a) we see that $\g$ is 
contained in the closure of $\g'$; since $\g\ne\g'$ it follows that $\dim\g'>\dim\g$, as required.

Note that (a),(b) provide, in the case where $C$ is elliptic, another characterization of $\Ph(C)$ which does
not rely on the partial order on $\uuG$.

\head 1. The main results\endhead
\subhead 1.1\endsubhead
In this section we assume that $p=2$.
Let $V$ be a $\kk$-vector space of finite dimension $\nn=2n\ge4$ with a fixed nondegenerate 
symplectic form $(,):V\T V@>>>\kk$ and a fixed quadratic form $Q:V@>>>\kk$ such that (i) or (ii) below holds:

(i) $Q=0$;

(ii) $Q\ne0$, $(x,y)=Q(x+y)-Q(x)-Q(y)$ for $x,y\in V$.
\nl
Let $Is(V)$ be the group consisting of all $g\in GL(V)$ such that
$(gx,gy)=(x,y)$ for all $x,y\in V$ and $Q(gx)=Q(x)$ for all $x\in V$ (a closed subgroup of $GL(V)$). 
Let $G$ be the identity component of $Is(V)$. Let $\cf$ be the set of all sequences 
$V_*=(0=V_0\sub V_1\sub V_2\sub\do\sub V_\nn=V)$ of subspaces of $V$ such that $\dim V_i=i$ for $i\in[0,\nn]$, 
$Q|_{V_i}=0$ and $V_i^\pe=V_{\nn-i}$ for all $i\in[0,n]$. 
Here, for any subspace $V'$ of $V$ we set $V'{}^\pe=\{x\in V;(x,V')=0\}$. 

\subhead 1.2\endsubhead
Let $p_1\ge p_2\ge\do\ge p_\s$ be a sequence in $\ZZ_{>0}$ such that $p_1+p_2+\do+p_\s=n$. 
(In the case where $Q\ne0$ we assume that $\s$ is even.)
For any $r\in[1,\s]$ we set $p_{\le r}=\sum_{r'\in[1,r]}p_{r'}$, $p_{<r}=\sum_{r'\in[1,r-1]}p_{r'}$. We fix 
$(V_*,V'_*)\in\cf\T\cf$ such that for any $r\in[1,\s]$ we have
$$\dim(V'_{p_{<r}+i}\cap V_{p_{<r}+i})=p_{<r}+i-r,\qua \dim(V'_{p_{<r}+i}\cap V_{p_{<r}+i+1})=p_{<r}+i-r+1\tag a$$
if $i\in[1,p_r-1]$;
$$\dim(V'_{p_{\le r}}\cap V_{\nn-p_{<r}-1})=p_{\le r}-r,\qua \dim(V'_{p_{\le r}}\cap V_{\nn-p_{<r}})
=p_{\le r}-r+1.\tag b$$
(Such $(V_*,V'_*)$ exists and is unique up to conjugation by $Is(V)$.)
A unipotent class $\g$ in $G$ is said to be adapted to $(V_*,V'_*)$ if for some $g\in\g$ we have $gV_i=V'_i$
for all $i$.

There is a unique unipotent conjugacy class $\g$ in $G$ such that $\g$ is adapted to $(V_*,V'_*)$ and some/any 
element of $\g$ has Jordan blocks of sizes $2p_1,2p_2,\do,2p_\s$. (The existence of $\g$ is proved in 
\cite{\WE, 2.6, 2.12}; the uniqueness follows from the proof of \cite{\WE, 4.6}.)

\proclaim{Theorem 1.3} Let $\g'$ be a unipotent conjugacy class in $G$ which is adapted to $(V_*,V'_*)$. Then 
$\g$ is contained in the closure of $\g'$ in $G$.
\endproclaim
The proof is given in 1.5-1.8 (when $Q=0$) and in 1.9 (when $Q\ne0$).

\subhead 1.4\endsubhead
Let $\ct$ be the set of sequences $c_*=(c_1\ge c_2\ge c_3\ge\do)$ in $\NN$ such that $c_m=0$ for $m\gg0$
and $c_1+c_2+\do=\nn$. For $c_*\in\ct$ we define $c_*^*=(c^*_1\ge c^*_2\ge c^*_3\ge\do)\in\ct$ by 
$c^*_i=|\{j\ge1;c_j\ge i\}|$ and we set $\mu_i(c_*)=|\{j\ge1;c_j=i\}|$ ($i\ge1$); thus we have
$$\mu_i(c_*)=c^*_i-c^*_{i+1}.\tag a$$
For $i,j\ge1$ we have
$$i\le c_j\text{ iff }j\le c^*_i.\tag b$$
For $c_*\in\ct$ and $i\ge1$ we have

(c) $\sum_{j\in[1,c^*_i]}(c_j-i)+\sum_{j\in[1,i]}c^*_j=\nn$.
\nl
Indeed the left hand side is  
$$\align&\sum_{j\ge1;i\le c_j}(c_j-i)+\sum_{j\in[1,i],k\ge 1;c_k\ge j}1
=\sum_{j\ge1;i\le c_j}(c_j-i)+\sum_{k\ge 1}\min(i,c_k)\\&
=\sum_{j\ge1;i\le c_j}(c_j-i)+\sum_{k\ge 1;i\le c_k}i+\sum_{k\ge 1;i>c_k}c_k\\&=
\sum_{j\ge1;i\le c_j}c_j+\sum_{k\ge 1;i>c_k}c_k=\sum_{j\ge1}c_j=\nn.\endalign$$
For $c_*,c'_*\in\ct$ and $i\ge1$:

(d) {\it we have $\sum_{j\in[1,i]}c^*_j=\sum_{j\in[1,i]}c'{}^*_j$ iff
$\sum_{j\in[1,c^*_i]}(c_j-i)=\sum_{j\in[1,c'{}^*_i]}(c'_j-i)$; 

we have $\sum_{j\in[1,i]}c^*_j\ge\sum_{j\in[1,i]}c'{}^*_j$ iff
$\sum_{j\in[1,c^*_i]}(c_j-i)\le\sum_{j\in[1,c'{}^*_i]}(c'_j-i)$.}
\nl
This follows from (c) and the analogous equality for $c'_*$.

For $c_*,c'_*\in\ct$ we say that $c_*\le c'_*$ if the following (equivalent) conditions are satisfied:

(i) $\sum_{j\in[1,i]}c_j\le\sum_{j\in[1,i]}c'_j$ for any $i\ge1$;

(ii) $\sum_{j\in[1,i]}c^*_j\ge\sum_{j\in[1,i]}c'{}^*_j$ for any $i\ge1$.
\nl
We show:

(e) {\it Let $c_*,c'_*\in\ct$ and $i\ge1$ be such that $c_*\le c'_*$,
$\sum_{j\in[1,i]}c^*_j=\sum_{j\in[1,i]}c'{}^*_j$ and $\mu_i(c_*)>0$. Then $c^*_i\le c'{}^*_i$ and 
$\mu_i(c'_*)>0$.}
\nl
We set $m=c^*_i,m'=c'{}^*_i$.
From $c_*\le c'_*$ we deduce $\sum_{j\in[1,i-1]}c^*_j\ge\sum_{j\in[1,i-1]}c'{}^*_j$ (if $i=1$ both sums are 
zero); using the equality $\sum_{j\in[1,i]}c^*_j=\sum_{j\in[1,i]}c'{}^*_j$ we deduce $c^*_i\le c'{}^*_i$ that is,
$m\le m'$. From (d) we have

$\sum_{j\in[1,m]}(c_j-i)=\sum_{j\in[1,m']}(c'_j-i)$.
\nl
Hence
$$\align&\sum_{j\in[1,m]}c_j=\sum_{j\in[1,m']}c'_j+(m-m')i\\&=\sum_{j\in[1,m]}c'_j+\sum_{j\in[m+1,m']}(c'_j-i)\ge
\sum_{j\in[1,m]}c'_j\ge\sum_{j\in[1,m]}c_j;\tag f\endalign$$
we have used $c_*\le c'_*$ and that for $j\in[m+1,m']$ we have $i\le c'_j$ (since $j\le c'{}^*_i$, see (b)).
It follows that the inequalities in (f) are equalities hence $c'_j=i$ for $j\in[m+1,m']$. Thus
$\mu_i(c'_*)\ge m-m'$. This completes the proof of (e) in the case where $m>m'$. Now assume that $m=m'$. From 
$c_*\le c'_*$ we have $\sum_{j\in[1,m-1]}c_j\le\sum_{j\in[1,m-1]}c'_j$. Using this and (d) we see that
$$\sum_{j\in[1,m]}(c_j-i)=\sum_{j\in[1,m]}(c'_j-i)\ge\sum_{j\in[1,m-1]}(c_j-i)+c'_m-i$$
hence $c_m-i\ge c'_m-i$. 
From $\mu_i(c_*)>0$ and $c^*_i=m$ we deduce that $c_m=i$. 
(Indeed by 1.4(b) we have $i\le c_m$; if $i<c_m$ then $i+1\le c_m$ and by 1.4(b) we have 
$m\le c^*_{i+1}\le c^*_i=m$ hence $c^*_{i+1}=c^*_i$ and $\mu_i(c_*)=0$, contradiction.)
Hence $c'_m\le i$.
Since $c'{}^*_i=m$ we have also $i\le c'_m$ (see (b)) hence $c'_m=i$. Thus $\mu_i(c'_*)>0$. This completes the 
proof of (e). 

\subhead 1.5\endsubhead
In this subsection (and until the end of 1.8)
we assume that $Q=0$. In this case we write $Sp(V)$ instead of $Is(V)=G$.
Let $u$ be a unipotent element of $Sp(V)$.
We associate to $u$ the sequence $c_*\in\ct$ whose nonzero terms are the sizes of the Jordan blocks of $u$.
We must have $\mu_i(c_*)=$even for any odd $i$. We also associate to $u$ a map 
$\e_u:\{i\in 2\NN;i\ne0,\mu_i(c_*)>0\}@>>>\{0,1\}$ as
follows: $\e_u(i)=0$ if $((u-1)^{i-1}(x),x)=0$ for all $x\in\ker(u-1)^i:V@>>>V$ and $\e_u(i)=1$ otherwise; we have
automatically $\e_u(i)=1$ if $\mu_i(c_*)$ is odd.
Now $u\m(c_*,\e_u)$ defines a bijection between the set of conjugacy classes of unipotent elements
in $Sp(V)$ and the set $\fS$ consisting of all pairs $(c_*,\e)$ where $c_*\in\ct$ is such that
$\mu_i(c_*)=$even for any odd $i$ and $\e:\{i\in2\NN;i\ne0,\mu_i(c_*)>0\}@>>>\{0,1\}$ is a function such that 
$\e(i)=1$ if $\mu_i(c_*)$ is odd. (See \cite{\SPA, I,2.6}). We denote by $\g_{c_*,\e}$ the unipotent class
corresponding to $(c_*,\e)\in\fS$.
For $(c_*,\e)\in\fS$ it will be convenient to extend $\e$ to a function $\ZZ_{>0}@>>>\{-1,0,1\}$ (denoted again 
by $\e$) by setting $\e(i)=-1$ if $i$ is odd or $\mu_i(c_*)=0$.

Now let $\g,\g'$ be as in 1.3. We write $\g=\g_{c_*,\e},\g'=\g_{c'_*,\e'}$ with $(c_*,\e),(c'_*,\e')\in\fS$. 
Let $g\in\g_{c'_*,\e'}$ be such that $gV_*=V'_*$ and let $N=g-1:V@>>>V$. To 
prove that $\g$ is contained in the closure of $\g'$ in $G$ it is enough to show that 

(a) $c_*\le c'_*$
\nl
and that for any $i\ge1$, (b),(c) below hold:

(b) $\sum_{j\in[1,i]}c^*_j-\max(\e(i),0)\ge\sum_{j\in[1,i]}c'{}^*_j-\max(\e'(i),0)$;

(c) if $\sum_{j\in[1,i]}c^*_j=\sum_{j\in[1,i]}c'{}^*_j$ and $c^*_{i+1}-c'{}^*_{i+1}$ is odd then $\e'(i)\ne0$.
\nl
(See \cite{\SPA, II,8.2}.)
From the definition we see that $c_i=2p_i$ for $i\in[1,\s]$, $c_i=0$ for $i>\s$ and from \cite{\WE, 4.6} we
see that $\e(i)=1$ for all $i\in\{2,4,6,\do\}$ such that $\mu_i(c_*)>0$. 

Now (a) follows from \cite{\WE, 3.5(a)}. Indeed in {\it loc.cit.} it is shown that for any $i\ge1$ we have 
$\dim N^iV\ge\L_i$ where 
$$\L_i=\sum_{j\ge1;i\le c_j}(c_j-i)=\sum_{j\in[1,c^*_i]}(c_j-i).$$
We have $\dim N^iV=\sum_{j\ge1;i\le c'_j}(c'_j-i)=\sum_{j\in[1,c'{}^*_i]}(c'_j-i)$ hence
by 1.4(d) the inequality $\dim N^iV\ge\L_i$ is the same as the inequality
$\sum_{j\in[1,i]}c^*_j\ge\sum_{j\in[1,i]}c'{}^*_j$. 

Note also that, by 1.4(d), 

(d) {\it we have $\sum_{j\in[1,i]}c^*_j=\sum_{j\in[1,i]}c'{}^*_j$ iff $\dim N^iV=\L_i$.}

\subhead 1.6\endsubhead
Let $i\ge1$. We show:

(a) {\it If $\mu_i(c_*)>0$ and $\sum_{j\in[1,i]}c^*_j=\sum_{j\in[1,i]}c'{}^*_j$ then $\e'(i)=1$.}
\nl
By 1.4(e) we have $\mu_i(c'_*)>0$. Since $\mu_i(c_*)>0$ we see that $i=2p_d$ for some $d\in[1,\s]$.
If $\mu_i(c'_*)$ is odd then $\e'(i)=1$ (by definition, since $i$ is even). Thus we may assume that 
$\mu_i(c'_*)\in\{2,4,6,\do\}$. From our assumption we have that $\dim N^iV=\L_i$ (see 1.5(d)).

Let $v_1,v_2,\do,v_\s$ be vectors in $V$ attached to $V_*,V'_*,g$ as in \cite{\WE, 3.3}.
For $r\in[1,\s]$ let $W_r,W'_r$ be as in \cite{\WE, 3.4}; we set $W_0=0$, $W'_0=V$.
From \cite{\WE, 3.5(b)} we see that $N^iW'_{d-1}=0$ at least if $d\ge2$; but the same
clearly holds if $d=1$. We have $v_d\in W'_{d-1}$ hence $N^{2p_d}v_d=0$ and
$$\align&(N^{2p_d-1}(v_d),v_d)=(N^{p_d}v_d,N^{p_d-1}v_d)=((g-1)^{p_d}v_d,(g-1)^{p_d-1}v_d)=\\&
(g^{p_d}v_d,v_d)=1.\endalign$$
(We have used that $(v_d,g^kv_d)=0$ for $k\in[-p_d+1,p_d-1]$ and $(v_d,g^{p_d}v_d)=1$, see \cite{\WE, 3.3(iii)}.)
Thus $\e'(i)=1$. This proves (a).

\subhead 1.7\endsubhead
We prove 1.5(b). It is enough to show that, if $\e(i)=1$ and $\e'(i)\le0$ then
$\sum_{j\in[1,i]}c^*_j\ge\sum_{j\in[1,i]}c'{}^*_j+1$. Assume this is not so. Then using 1.5(a) we have
$\sum_{j\in[1,i]}c^*_j=\sum_{j\in[1,i]}c'{}^*_j$. Since $\e(i)=1$ we have $\mu_i(c_*)>0$; using 1.6(a) we see that
$\e'(i)=1$, a contradiction. Thus 1.5(b) holds.

\subhead 1.8\endsubhead
We prove 1.5(c). If $i$ is odd then $\e'(i)=-1$, as required.
Thus we may assume that $i$ is even. Using 1.5(a) and 1.4(e) we see that $c^*_i\le c'{}^*_i$.

Assume first that $c^*_i=c'{}^*_i$. From $\mu_i(c_*)=c^*_i-c^*_{i+1}$, $\mu_i(c'_*)=c'{}^*_i-c'{}^*_{i+1}$ we
deduce that $\mu_i(c_*)-\mu_i(c'_*)=c'{}^*_{i+1}-c^*_{i+1}$ is odd.
If $\mu_i(c'_*)$ is odd we have $\e'(i)=1$ (since $i$ is even); thus we have
$\e'(i)\ne0$, as required. If $\mu_i(c'_*)=0$ we have $\e'(i)=-1$; thus we have $\e'(i)\ne0$, as required. 
If $\mu_i(c'_*)\in\{2,4,6,\do\}$ then $\mu_i(c_*)$ is odd so that $\mu_i(c_*)>0$ and then 1.6(a) shows that
$\e'(i)=1$; thus we have $\e'(i)\ne0$, as required. 

Assume next that $c^*_i<c'{}^*_i$. By 1.5(a) we have $\sum_{j\in[1,i+1]}c^*_j\ge\sum_{j\in[1,i+1]}c'{}^*_j$; 
using the assumption of 1.5(c) we deduce that $c^*_{i+1}\ge c'{}^*_{i+1}$. Combining this with $c^*_i<c'{}^*_i$ 
we 
deduce $c^*_i-c^*_{i+1}<c'{}^*_i-c'{}^*_{i+1}$ that is, $\mu_i(c_*)<\mu_i(c'_*)$. It follows that $\mu_i(c'_*)>0$.
If $\mu_i(c_*)>0$ then by 1.6(a) we have $\e'(i)=1$; thus we have $\e'(i)\ne0$, as required. Thus we can assume 
that $\mu_i(c_*)=0$. We then have $c^*_i=c^*_{i+1}$ and we set $\d=c^*_i=c^*_{i+1}$. As we have seen earlier, we 
have $c^*_{i+1}\ge c'{}^*_{i+1}$; using this and the assumption of 1.5(c) 
we see that $c^*_{i+1}-c'{}^*_{i+1}=2a+1$ where 
$a\in\NN$. It follows that $c'{}^*_{i+1}=\d-(2a+1)$. In particular we have $\d\ge2a+1>0$.

If $k\in[0,2a]$ we have $c'_{\d-k}=i$. (Indeed, assume that $i+1\le c'_{\d-k}$; then by 1.4(b) we have 
$\d-k\le c'{}^*_{i+1}=\d-(2a+1)$ hence $k\ge 2a+1$, a contradiction. Thus $c'_{\d-k}\le i$. On the other hand,
$\d=c^*_i<c'{}^*_i$ implies by 1.4(b) that $i\le c'_\d$. Thus $c'_{\d-k}\le i\le c'_\d\le c'_{\d-k}$ hence
$c'_{\d-k}=i$.)

Using 1.4(b) and $c'{}^*_{i+1}=\d-(2a+1)$ we see that $c'_{\d-(2a+1)}\ge i+1$ (assuming that $\d-(2a+1)>0$). Thus
the sequence $c'_1,c'_2,\do,c'_\d$ contains exactly $2a+1$ terms equal to $i$, namely
$c'_{\d-2a},\do,c'_{\d-1},c'_\d$.

We have $i>c_{\d+1}$. (If $i\le c_{\d+1}$ then from 1.4(b) we would get $\d+1\le c^*_i=\d$, a contradiction.)

Since $\d>0$, from $c^*_i=\d$ we deduce that $i\le c_\d$ (see 1.4(b)); since 
$\mu_i(c_*)=0$ we have $c_\d\ne i$ hence $c_\d>i$. From the assumption of 1.5(c) we see that $\dim N^iV=\L_i$ 
(see 1.5(d)). Using this and $c_\d>i>c_{\d+1}$ we see that \cite{\WE, 3.5} is applicable and gives that 
$V=W_\d\op W_\d^\pe$ and $W_\d$, $W_\d^\pe$ are $g$-stable; moreover, $g:W_\d@>>>W_\d$ has exactly $\d$ Jordan 
blocks and each one has size $\ge i$ and $g:W_\d^\pe@>>>W_\d^\pe$ has only Jordan blocks of size $\le i$. Since 
the $\d$ largest numbers in the sequence $c'_1,c'_2,\do$ are $c'_1,c'_2,\do,c'_\d$ we see that the sizes of the 
Jordan blocks of $g:W_\d@>>>W_\d$ are $c'_1,c'_2,\do,c'_\d$. Since the last sequence contains an odd number 
($=2a+1$) of terms equal to $i$ we see that $\e_{g|_{W_d}}(i)=1$. (Note that $(,)$ is a nondegenerate symplectic 
form on $W_d$ hence $\e_{g|_{W_d}}(i)$ is defined as in 1.5.) Hence there exists $z\in W_d$ such that $N^iz=0$ 
and $(z,N^{i-1}z)=1$. This implies that $\e_g(i)=1$ that is $\e'(i)=1$. This completes the proof of 1.5(c) and 
also completes the proof of Theorem 1.3 when $Q=0$.

\subhead 1.9\endsubhead
In this subsection we assume that $Q\ne0$. Let $\g,\g'$ be as in 1.3. We denote by $\g_1,\g'_1$ the 
$Is(V)$-conjugacy class containing $\g,\g'$ respectively; let $\g_2,\g'_2$ be the $Sp(V)$-conjugacy class 
containing $\g_1,\g'_1$ respectively. Note that Theorem 1.3 is applicable to $\g_2,\g'_2$ instead of $\g,\g'$ 
and with $G$ replaced by the larger group $Sp(V)$. Thus we have that $\g_2$ is contained in the closure of 
$\g'_2$ in $Sp(V)$ and then, using \cite{\SPA, II,8.2}, we see that $\g_1$ is contained in the closure of $\g'_1$
in $Is(V)$. We have $\g_1=\g$ (see \cite{\SPA, I,2.6}). If $\g'_1=\g'$ it follows that $\g$ is contained in the 
closure of $\g'$ in $G$, as required. If $\g'_1\ne\g'$ then $\g'_1=\g'\sqc\g''$ where $\g''=r\g'r\i$ ($r$ is a 
fixed element in $Is(V)-G$). We see that either $\g$ is contained in the closure of $\g'$ or in the closure of 
$r\g'r\i$. In the last case we have that $r\i\g r$ is contained in the closure of $\g'$. But $r\i\g r=\g$ so
that in any case $\g$ is contained in the closure of $\g'$. This completes the proof of Theorem 1.3 when $Q\ne0$.

\enddocument